\documentclass{article}
\usepackage[utf8]{inputenc}
\usepackage[colorlinks=true,allcolors=blue]{hyperref}
\usepackage{geometry,graphicx}
\usepackage{amsmath, amssymb, amsthm}
\usepackage{booktabs}
\usepackage{threeparttable}

\usepackage{graphicx} 
\usepackage{amsmath,amssymb,amsthm}
\usepackage{algpseudocode}
\usepackage{algorithm}
\usepackage{fullpage}
\usepackage{hyperref}
\usepackage{xcolor}
\usepackage{mleftright}

\newcommand{\R}{{\mathbb{R}}}

\newcommand{\eps}{\varepsilon}

\newcommand{\inr}[2]{\left\langle#1,#2\right\rangle}

\newcommand{\wrt}{\,\textnormal d}
\newcommand{\deriv}[2]{\frac{\wrt^{#2}}{\wrt #1^{#2}}}

\newcommand{\abs}[1]{\left|#1\right|}
\newcommand{\magn}[1]{\left\|#1\right\|}

\newcommand{\pare}[1]{\mleft(#1\mright)}

\newcommand{\set}[1]{{\left\{{#1}\right\}}}

\newcommand{\ceil}[1]{\left\lceil#1\right\rceil}
\newcommand{\bmat}[1]{\begin{bmatrix}#1\end{bmatrix}}

\newcommand{\indicate}[1]{{\mathbf1_{#1}}}

\DeclareMathOperator*{\E}{\mathbb{E}}

\newtheorem{theorem}{Theorem}[section]
\newtheorem{lemma}[theorem]{Lemma}

\theoremstyle{definition}
\newtheorem{definition}{Definition}[section]

\theoremstyle{remark}
\newtheorem*{remark}{Remark}

\algdef{SE}[REPEATN]{RepeatN}{End}[1]{\algorithmicrepeat\ #1 \textbf{times}}{\algorithmicend}

\theoremstyle{definition}

\newcommand{\rvec}{{\mathbf r}}

\newcommand{\svec}{{\mathbf s}}
\newcommand{\e}{{\mathbf e}}
\newcommand{\len}{\textnormal{len}}

\DeclareMathOperator{\sigg}{\sigma}
\newcommand{\sig}[1]{\sigg\mleft(#1\mright)}

\title{Achieving the Highest Possible Elo Rating}
\author{Rikhav Shah}
\date{February 2024}

\begin{document}

\maketitle

\begin{abstract}
Elo rating systems measure the approximate skill of each competitor in a game or sport. A competitor's rating increases when they win and decreases when they lose. Increasing one's rating can be difficult work; one must hone their skills and consistently beat the competition.   Alternatively, with enough money you can rig the outcome of games to boost your rating. This paper poses a natural question for Elo rating systems: say you manage to get together $n$ people (including yourself) and acquire enough money to rig $k$ games.  How high can you get your rating, asymptotically in $k$?  In this setting, the people you gathered aren't very interested in the game, and will only play if you pay them to. This paper resolves the question for $n=2$ up to constant additive error, and provide close upper and lower bounds for all other $n$, including for $n$ growing arbitrarily with $k$. There is a phase transition at $n=k^{1/3}$: there is a huge increase in the highest possible Elo rating from $n=2$ to $n=k^{1/3}$, but (depending on the particular Elo system used) little-to-no increase for any higher $n$. Past the transition point $n>k^{1/3}$, the highest possible Elo is at least $\Theta(k^{1/3})$. The corresponding upper bound depends on the particular system used, but for the standard Elo system, is $\Theta(k^{1/3}\log(k)^{1/3})$.

\end{abstract}

\section{Introduction}
The Elo Rating system was proposed by Arpad Elo in the mid 20th century to estimate the relative skill of chess players \cite{original_elo}.  It was quickly adopted by the international chess community, and in the decades since has seen adoption in many competitive contexts.
This paper considers a simple combinatorial question about the Elo system.
To the author's knowledge, this is the first time this question has been posed in print: \textbf{given $n$ players starting with equal rating, what is the highest a player could be rated after a total of $k$ games are played?}

We begin with a definition of the system then provide its motivation.
Each player is given some `rating' value (measured in `points' or simply `Elo'), which updates as they play games.  These rating points are somewhat analogous to poker chips: when player $A$ and player $B$ play a game, they each place some of their rating points into a pot.  In the case of a draw, the players split the pot evenly.  If one player wins, they take the entire pot.
The heart of the Elo system is dictating how many points each player must ante up.
To do so, each implementation of the system specifies a `pot function' $\sigma$ satisfying
\begin{enumerate}
    \item $\sigma$ is non-negative and monotonically increasing, and
\vspace{-.5cm}
\begin{equation}\label{list:pot_function_requirements}\end{equation}
\vspace{-1.0cm}
    \item $\sigma(z)+\sigma(-z)=1$ for all $z\in\R$.
\end{enumerate}
Let $r_A$ and $r_B$ be the ratings of players $A$ and $B$ respectively. 
When players $A$ and $B$ play, the number of points they ante up are $K\cdot\sigma(r_A-r_B)$ and $K\cdot\sigma(r_B-r_A)$ respectively, for a total pot size of $K$. Players are allowed to go into debt if they don't have the required points, i.e. negative ratings are perfectly fine.
The value of $K$ itself is a parameter of the system. 
The resulting rating updates for different outcomes of a game between $A$ and $B$ are as follows:
\begin{equation}
\label{eq:update_rule}
\begin{tabular}{ l|c|c|c| }
& $A$ wins & $A$ and $B$ draw & $B$ wins \\
\hline
\rule{0pt}{4ex}
$r_A'=$ & $r_A+K\cdot\sigma(r_B-r_A)$ & $r_A+\frac K2\cdot\pare{\sigma(r_B-r_A)-\sigma(r_A-r_B)}$ & $r_A-K\cdot\sigma(r_A-r_B)$
\\
[2ex]
\hline
\rule{0pt}{4ex}
$r_B'=$ & $r_B-K\cdot\sigma(r_B-r_A)$ & $r_B+\frac K2\cdot\pare{\sigma(r_A-r_B)-\sigma(r_B-r_A)}$ & $r_B+K\cdot\sigma(r_A-r_B)$
\\
[2ex]
\hline
\end{tabular}
\end{equation}
The motivation for this system comes from thinking of the outcome of a game as a random variable. For some symmetric random variable $\eta$, the event $r_A-r_B+\eta>0$ is recorded as a victory for $A$, $r_A-r_B+\eta<0$ as a victory for $B$, and $r_A-r_B+\eta=0$ as a draw.
Set
\begin{equation}
    \label{eq:prob_interp}
    \sigma(z)=\Pr(\eta<z)+\frac12\Pr(\eta=z)
\end{equation}
so that $\sigma$ is a kind of `symmetrized' cumulative distribution function of $\eta$ (which coincides with the usual cumulative distribution function when $\eta$ has no atoms).
Given $\eta$, setting $\sigma$ to satisfy (\ref{eq:prob_interp}) guarantees it satisfies (\ref{list:pot_function_requirements}). Conversely, given any $\sigma$ satisfying (\ref{list:pot_function_requirements}), one can define a symmetric random variable $\eta$ satisfying (\ref{eq:prob_interp}), if some of the mass of $\eta$ is allowed to be at infinity (i.e. $\lim_{z\to\infty}\Pr(\eta\le z)$ need not be $1$ and $\lim_{z\to-\infty}\Pr(\eta<z)$ need not be $0$).
Under this probabilistic model, observe that $\E(r'_A)=r_A$ and $\E(r'_B)=r_B$. So one gains points for performing `better than expected' and loses points for performing `worse than expected'.
One additional natural assumption is that $\eta$ has a finite expectation, though it is not required.
Given some real-world game, one should attempt to pick $\eta$ so that these estimated probabilities match the empirical win-loss rates observed.
The proposal by Elo in the 1960s was to take $\eta$ to be Gaussian, citing the ubiquity of the normal distribution in nature \cite{suggest_gaussian}. However, the community soon decided a logistic random variable was more suitable, leading to the pot function of $\sigma(z)=\frac1{1+e^{-cz}}$ for some constant $c$ \cite{original_elo}. The International Chess Federation (FIDE) has long used $c=\log(10)/400\approx 5.76\times10^{-3}$ \cite[Chapter~B02]{fide_handbook}, though recent analysis suggests $c=\frac56\cdot\log(10)/400\approx4.80\times10^{-3}$ reflects real-world chess data much better \cite{adjusting_c}.
Depending on various factors, FIDE uses $K$ in the range of $10$ to $40$. The range of ratings exhibited by human players is roughly $0$ to $3000$ \cite{fide_handbook,fide_ratings}.
This paper considers a generic pot function $\sigma$, and applies the results to several specific families of pot functions listed in Table \ref{tab:results}.


There are additional complications in real-world implementations of Elo. For legibility and practicality, fractional and negative rating points are avoided by scaling and shifting points up and rounding to the nearest integer, and by imposing an artificial floor on possible ratings (by gifting a player points if they would otherwise dip below the floor).
The total size of the pot $K$ may also vary depending on various factors, such as how many games each player has played before. For example, $K$ may be large for a new player to facilitate faster convergence of their rating to their true skill level, and may decrease over time to reduce arbitrary fluctuations for experienced players. Sometimes rating updates are batched. That is, one accumulates their pot winnings and losses over several games, and updates their rating once at the end. Often times all the games played at a single tournament are batched together in that manner. All these details and even more complications are outlined thoroughly by different organizations implementing Elo; one may read about them in, for example, the FIDE Handbook \cite{fide_handbook}.


\subsection{Setting and results}
This paper considers $n$ players starting with equal ratings. %
Players' ratings update on a game-by-game basis according to (\ref{eq:update_rule}) with $K=1$ kept fixed.
In particular, both fractional and negative points are possible and rating updates are not batched.
Since the dynamics depend only on the difference between ratings, we can take everyone's initial rating to be 0 without loss of generality.
Asymptotically in $k$, we seek the highest one of the $n$ players may be rated after a total of $k$ games are played amongst all of them. This question is interesting both for fixed $n$, and for $n$ allowed to grow with $k$.
We find a phase transition at $n=k^{1/3}$: there is a huge increase in the highest possible Elo rating from $n=2$ to $n=k^{1/3}$, but (depending on $\sigma$) little-to-no increase for any higher $n$.
This paper is organized into four main sections.
\begin{enumerate}
    \item[] Section \ref{sec:nis2}: the highest Elo problem for $n=2$ is resolved up to constant additive error.
    \item[] Section \ref{sec:opt_strat_large_n}: a lower bound is provided for each $n$ by finding a family of strategies which achieve a highest rating of $\Theta\pare{\min(n,k^{1/3})}$ for any pot function. Note the bound does not improve past $n=\Omega\pare{k^{1/3}}$.
    \item[] Section \ref{sec:opt}: an upper bound on the highest possible rating for each $n$ is provided with a mild natural assumption on $\sigma$. 
    \item[] Section \ref{sec:dis}: open questions and a surprising connection to the maximum overhang problem are discussed.
\end{enumerate}
The gap between the upper and lower bounds depends on which $\sigma$ is used. In particular, the key quantity is the left-tail behavior of $\sigma$.
A quickly decaying tail corresponds both to a small asymptotic rate in the case of $n=2$, and to nearly matching upper and lower bounds for general $n$.
Our main result is stated in terms of the following function, which quantifies the rate of decay of the left-tail of $\sigma$:
\begin{equation}
\label{eq:f_declr}
    f(x)=\int_0^{x}\frac1{\sigma(-\tau)}\,\wrt\tau
\end{equation}
Note that faster the left tail of $\sigma$ decays to $0$, the faster $f$ diverges to infinity.
Note that non-negativity and monotonicity of $\sigma$ means that $f$ is increasing and convex.
This main theorem summarizes the results from each section.
\begin{theorem}
\label{THM:FINAL}
    Let $R(n,k)$ be the highest possible Elo rating achievable with $k$ games $n$ players starting with $0$ Elo points.
    Fix any pot function $\sigma$ and let $f$ be as in (\ref{eq:f_declr}).
    Then 
    \[\frac12f^{-1}(2k)\le R(2,k)\le\frac12f^{-1}(2k-2)+\frac52.\]
    Now assume that $\sup_zz\,\sig{-z}<\infty$.
    Fix $n=h(k)$.
    If $h(k)=\Omega(k^{1/3})$, then there exists constants $C_1,C_2$ such that for sufficiently large $k$, one has
    \[C_1 k^{1/3}\le R(h(k), k) \le C_2 k^{1/3}f^{-1}(k)^{1/3}.\]
    Furthermore, if $h(k)=o(k^{1/3})$ then
    \[C_1 h(k)\le R(h(k), k) \le C_2 h(k) f^{-1}(k/h(k)).\]
\end{theorem}
\begin{proof}
    $R(2,k)$ handled in Section \ref{sec:nis2}. The lower bounds for $R(h(k),k)$ are in Section \ref{sec:opt_strat_large_n} and the upper bounds are in Section \ref{sec:opt}.
\end{proof}
\begin{remark}
    The results hold even if the players all start with distinct ratings, provided that all the initial ratings fall into some bounded interval independent of $n$ and $k$, and that there's some point around which the inital ratings are symmetric. That is, there exists some value $r_0$ such that for every player rated $r$ there's one rated $2r_0-r$.
\end{remark}

Table \ref{tab:results} lists the value of $f^{-1}$ for several natural families of $\sigma$ and names the corresponding $\eta$ where appropriate.
For the logistic pot function, $f^{-1}$ is logarithmic and Theorem \ref{THM:FINAL} implies
\[R(2,k)=\frac1{2c}\log(2k)+O(1)\quad\text{and}\quad
R(n,k)=\widetilde\Theta\pare{\min\pare{n,k^{1/3}}}
\]
where $\widetilde\Theta$ suppresses log factors.
For any pot function that eventually hits 0 (for example, the `uniform' pot function $\max(0,\min(1,cz+1/2))$), let $x$ be such that $\sig{-x}=0$ and note $f$ has a vertical asymptote at $x$. This means $f^{-1}$ is bounded, so Theorem \ref{THM:FINAL} implies
\[R(n,k)=\Theta\pare{\min\pare{n,k^{1/3}}}\]
where $\widetilde\Theta$ only suppresses constant factors.
At the other extreme, $\sigma$ need not converge to 0 at all. In this case, $f^{-1}(k)=\Theta(k)$. The upper bounds in Theorem \ref{THM:FINAL} no longer hold, but just by noticing the total pot size for each game is bounded, one sees $R(n,k)=\Theta(f^{-1}(k))=\Theta(k)$ anyway.




\begin{table}[h]\centering
\begin{tabular}{@{}lll@{}}
\toprule
$\eta$ & $\sigma(z)$ & $f^{-1}(k)$ \\
\midrule
\midrule
Logistic & $\frac1{1+e^{-cz}}$ & $\frac1{c}\log(ck)\cdot(1+O(1/k))$ \\
\midrule
--& $\frac12\frac{cz}{(1+|cz|^p)^{1/p}}+\frac12$ & $\pare{\frac1{2c^p}\cdot\frac{p+1}{p}}^{\frac1{p+1}}\cdot k^{\frac1{p+1}}\cdot(1+o(1))$ \\
\midrule
Gaussian & $\frac12\,\text{erf}(cz/\sqrt 2)+\frac12$ & $\frac1{c}\sqrt{\log k}\cdot(1+o(1))$ \\
\midrule
Uniform & $\min\pare{1,\max\pare{0,cz+\frac12}}$ & $\frac1{2c}-\frac{e^{-ck}}{2c}$ \\
\midrule
Cauchy & $\frac1\pi\arctan(cz)+\frac12$ & $\sqrt{2/c\pi}\cdot k^{1/2}+O(1)$ \\
\midrule
-- & $\frac c2\,\text{sign}(x)+\frac12$ & $\frac2{1-c}\cdot k$\\   
\bottomrule
\end{tabular}
\caption{
The value of $f^{-1}$ for some selected families of pot functions. Most of the families are parameterized by some constant $c$, which correlates with the slope of each $\sigma$ at 0.
The logistic pot function in the top row is the usual one used by real-world implementations of Elo.
}
\label{tab:results}
\end{table}

\section{Case of $n=2$}
\label{sec:nis2}
When there are only two players, $H$ and $L$, as their ratings $r_H$ and $r_L$ grow further apart, the fewer rating points the higher player $H$ can earn from the lower $L$ each time.
If $\sigma(-x)=0$ for some $x$, the difference in the ratings of $H$ and $L$ cannot exceed $x+2$. To see this, observe that a game can only increase the value of $r_H-r_L$ if $|r_H-r_L|<x$, and furthermore, $r_H$ and $r_L$ can each change by at most one point in each game. Since $r_H+r_L=0$ is invariant, this immediately implies
\[r_H\le\frac x2+1\]
independently of how many games are played.
On the other hand, if $\sigma(-z)>0$ for all $z$, then $r_H-r_L$ will diverge. To see this first note monotonicity of $\sigma$ means $\sigma$ is bounded away from 0 on each compact interval. Then, note $H$ beating $L$ will change the value of $r_H-r_L$ by applying the map $g:z\mapsto z+2\sigma(-z)$. The orbit can only converge if $\sigma(-g^{j}(z))$ converges to 0, which cannot happen for bounded $z$ by assumption.
In the former case, we may still ask how quickly the bound $r_H\le\frac z2+1$ is achieved, and in the latter case we ask how quickly the orbit diverges.
Intuitively, the faster $\sigma(-z)$ decays to 0 as $z\to\infty$, the slower the rate. Indeed, as the next theorem makes explicit, there is a simple expression for the rate in terms of the $f$ defined in (\ref{eq:f_declr}). Note $f$ is only well defined for $\set{x:\sigma(-x)>0}$. Nevertheless,
it's continuous and strictly increasing, and its range is all of $\R$. This implies for any $\sigma$ that it has a well defined inverse $f^{-1}$ defined on all of $\R$. The following theorem takes advantage of that fact to unify the analysis for every $\sigma$.


\begin{theorem}
\label{thm:nis2}
Fix any pot function $\sigma$. Let $r(t)$ be the highest possible rating after $t$ games with two players. Then
\[\frac{f^{-1}(2t)}2\le r(t)\le \frac{f^{-1}(2t-2)}2+\frac52.\]
\end{theorem}
\begin{proof}
Since the ratings of the two players sum to zero, it suffices to keep track of just the higher rated player.
Call the players $H$ and $L$ for `higher' and `lower'.
Let $r(t)$ be the rating of $H$ after $t$ games. Then $-r(t)$ is the rating of $L$ after $t$ games.
Under the assumption that $H$ wins every game, we have the simple recurrence
\begin{align*}
r(0)&=0 &  r(t+1)&=r(t)+\sigma(-r(t)-r(t))\\&&&=r(t)+\sigma(-2r(t)).
\end{align*}
This can be viewed as running Euler's method on the differential equation
\[y(0)=0\quad\quad\quad y'=\sigma(-2y)\]
with a step size of $1$. Since $y'$ is positive and monotonically decreasing in $y$, the Euler approximation upper bounds the exact solution. That is,
\[r(t)\ge y(t).\]
Using separate and integrate, one sees the exact solution to the differential equation is\[y(t)=\frac12f^{-1}(2t).\]
This establishes the lower bound. For the upper bound, we cannot assume that the optimal strategy is for $H$ to win every game. In particular, the function $z\mapsto z+\sigma(-2z)$ need not be monotone, so we cannot exclude the possibility one can achieve a higher rating by first losing a game and `slingshotting' to a higher rating exploiting the fact that $\sigma(-2y)$ is larger for smaller $y$.
Construct a sequence $x_j$ based on the recurrence
\begin{equation}
\label{eq:ub_recurrence}
x_0=1\quad\quad\quad x_{j+1}=x_j+\sigma(-2x_{j}+2).
\end{equation}
This sequence defines a partition of the positive real line,
\((0,x_0],(x_0,x_1],(x_1,x_2],(x_2,x_3]\cdots.\)
Again we take $r(t)$ to be the rating of $H$ after $t$ games, but instead of a strict recurrence we only have an inequality $r(t+1)\le r(t)+\sigma(-2r(t))$, which is tight only when $H$ wins game $t+1$. We now relate $r(t)$ to $x_j$.
If $r(t)\le0$, then $r(t+1)\le r(t)+1\le1=x_0$. If $r(t)\in(0,x_0]$, then $r(t+1)\le r(t)+\sigma(-2r(t))\le x_0+\sigma(-2x_0+2)=x_1$.
If $r(t)\in(x_{j-1},x_j]$ for some $j$, then
\begin{align*}
    r(t+1)
  &\le r(t)+\sigma(-2r(t))
\\&\le x_j+\sigma(-2x_{j-1})
\\&\le x_j+\sigma(-2x_{j}+2)
\\&=x_{j+1}.
\end{align*}
In each case, we see that
\[r(t)\le x_j\implies r(t+1)\le x_{j+1}.\]
Since $r(0)=0<1=x_0$, we have $r(t)\le x_t$ regardless of the sequence of wins and losses.
Now similar to before, $x_t$ can be viewed as the result of running Euler's method with a step size of $1$ on the differential equation
\begin{equation}
\label{eq:diffeq2}
    y'=\sigma(-2y+2).
\end{equation}
Using separate and integrate, see that the solution to the differential equation is
\[y(t)=\frac12f^{-1}(2t)+1.\]
Again the Euler approximation is an upper bound to the exact solution. That is,
\begin{equation}
\label{eq:xt_lb}
x_t\ge y(t).
\end{equation}
Suprisingly, this lower bound leads us to an upper bound. We first plug the recurrence (\ref{eq:ub_recurrence}) back into itself to express $x_t$ as a sum, then use (\ref{eq:xt_lb}) with monotonicity of $\sigma$.
\begin{equation}
\label{eq:xt_ub}
x_t=x_0+\sum_{j=0}^{t-1}\sigma(-2x_j+2)\le x_0+\sum_{j=0}^{t-1}\sigma(-2y(j)+2).
\end{equation}
The summand is the same as the right hand side of the differential equation (\ref{eq:diffeq2}). Then the sum is Reimann sum of step size $1$. But since the summand is monotone, we can bound the Reimann sum by the integral.
\begin{equation}
\label{eq:reimansum}
\sum_{j=0}^{t-1}\sigma(-2y(j)+2)
=\sum_{j=0}^{t-1}y'(j)
\le y'(0)+\int_{0}^{t-1}y'(s)\wrt s
=y'(0)+y(t-1).
\end{equation}
Note $y'(0)=1/2$ and $x_0=1$, so combining (\ref{eq:xt_ub}) and (\ref{eq:reimansum}) gives the final upper bound of
\[x_t\le1+\frac12+y(t-1)=1+\frac12+1+\frac12f^{-1}(2t-2).\]

\end{proof}

\section{Lower bound for general $n$}
\label{sec:opt_strat_large_n}
This section describes two strategies for any number of players and games.
The first strategy does not depend at all on the pot function used, but requires that everyone's initial ratings be exactly equal.
The second strategy has a small dependence on the pot function used, but works for any symmetric list of initial ratings.
Both strategies produce higher ratings for larger $n$, up to $n=\Theta(k^{1/3})$, at which point the asymtotic highest rating in $k$ is $\Theta(k^{1/3})$.
The first strategy is very simple: \textbf{pick any pair of players with equal rating and have one beat the other. Repeat until all players have a distinct rating or $k$ games have been played.}
This strategy is guaranteed to produce a player of either very high or very low rating. If it produces a player of very low rating, simply re-do the strategy picking the same sequence of pairs of players but have the opposite player win. Since game outcomes are symmetric, this will produce a player of high rating instead.

\begin{theorem}
\label{thm:construction}
The aforementioned strategy achieves a highest rating of
\[\min\pare{\frac{k^{1/3}}2,\frac n4}\]
when all players start with $0$ rating points.
In particular, the rating achieved is $\Omega(k^{1/3})$ for $n=\Omega(k^{1/3})$.
\end{theorem}
\begin{proof}
Let $\rvec(j)$ be the multiset of player ratings after $j$ games are played and $\abs{\rvec(j)}$ its element absolute value. We claim that every entry of $\rvec(j)$ is a half-integer. To see this, note $\sigma(0)=1/2$ for every pot function so player ratings change in increments of $1/2$ and all start at $0$.
We next claim
\begin{equation}
\label{eq:r_squared_tracks_k}
\sum_{r\in\rvec(j)}r^2=j/2.
\end{equation}
We argue inductively. The claim is clearly true initially when all ratings are 0. When two players rated $a$ points play, the increase in the functional is
\begin{align*}
\Delta=(a+1/2)^2+(a-1/2)^2-a^2-a^2=1/2.
\end{align*}
Suppose the strategy terminates when $k$ games have been played.
Then
\[\frac k2=\sum_{r\in\rvec(k)}r^2\le n\max\abs{\rvec(k)}^2
\implies \max\abs{\rvec(k)}\ge\sqrt{k/2n}.\]
Now suppose the strategy terminates when all players have distinct ratings.
Then since all ratings are half integers, that means $\max\abs{\rvec(j)}\ge n/4$.
Without a guarantee on which way the strategy terminates, we get the worst of both bounds. So the maximum rating is at least $\min\pare{\sqrt{k/2n},n/4}$. If $n\le2k^{1/3}$, the minimum is $n/4$. If $n>2k^{1/3}$, we can rerun the strategy ignoring all but the first $2k^{1/3}$ players, giving a bound of exactly $k^{1/3}/2$.
\end{proof}
\begin{remark}
Theorem \ref{thm:construction} is unfortunately very brittle. If the initial ratings are perturbed slightly, we can no longer assume all ratings are half-integers, nor can we often expect players to have equal rating. For
\[\sigma(x)=\frac12\,\textnormal{sign}(x)+\frac12=\begin{cases}1&x>0\\1/2&x=0\\0&x<0\end{cases},\]
this is indeed the best we can do. If all players have slightly different ratings, no Elo can be transferred between them. However, there's a robust version of the strategy that works for every other $\sigma$: first fix some $\delta>0$ for which $\sigma(-\delta)>0$.
\textbf{Pick any pair of players whose ratings are within $\delta$ of each other and have the higher rated player beat the lower rated player. Repeat until no two players are within $\delta$ rating points or until $k$ games have been played.} Note for $\delta=0$, one recovers the original strategy. As before, we may have to flip everyone's rating to ensure we end with a very high rating, as opposed to a very low rating. We need not assume that the initial ratings are all 0 anymore, but we do need to allow ourselves this possible reflection.
\end{remark}
{
\begin{theorem}
\label{thm:robust_construction}
Fix any $\vec r=(r_1,\cdots,r_n)$. The aforementioned strategy achieves a highest rating of
\[
\min\pare{\delta^{1/3}\sigma(-\delta)^{1/3}\cdot k^{1/3},\,\frac\delta2\cdot n}
\]
for at least one choice of the initial ratings $\rvec(0)=\vec r$ or $\rvec(0)=-\vec r$.
\end{theorem}
\begin{proof}
    The proof is very similar to the proof of Theorem \ref{thm:construction}. Again we consider the functional $\sum_{r\in\rvec(j)}r^2$. In a game between players rated $a<b$, the change in the functional is
\begin{align*}
\Delta
  &=(b+\sigma(a-b))^2+(a-\sigma(a-b))^2-b^2-a^2
\\&=2b\sigma(a-b)-2a\sigma(a-b)+2\sigma(a-b)^2
\\&=2(b-a)\sigma(a-b)+2\sigma(a-b)^2.
\\&\ge2\sig{-\delta}^2.
\end{align*}
Suppose the strategy terminates when $k$ games have been played.
Then
\[2\sig{-\delta}\cdot k\le\sum_{r\in\rvec(k)}r^2\le n\max\abs{\rvec(k)}^2
\implies \max\abs{\rvec(k)}\ge\sqrt{2\sig{-\delta}k/n}.\]
Now suppose the strategy terminates when no two players are within $\delta$ rating points. Since there are a total of $n$ players, that means $\max\abs{\rvec(j)}\ge\delta n/2$.
As before, without a guarantee on which way the strategy terminates, we get the worst of both bounds. When
\[n\le {\frac{2\sig{-\delta}^{1/3}}{\delta^{2/3}}}\cdot k^{1/3},\]
the smaller bound is $\delta n/2$. When $n$ is larger, one simply ignores the excess players, as before.
\end{proof}
\begin{remark}
The theorem is best applied with the vector of initial ratings is symmetric, i.e. $\vec r$ is a permutation of $-\vec r$. In fact, by shifting all the points up or down by a constant amount, they can be symmetric around any (constant valued) point, i.e. $\vec r$ is a permutation of $(r_0,\cdots,r_0)-\vec r$.
\end{remark}
\begin{remark}
This bound is stronger the heavier the tail of $\sigma$. Consider, for instance
\(
\sigma(z)=\frac{1}{1+e^{-cz}}
\)
as $c\to\infty$, or any other family of $\sigma$ approaching the pathological $\frac12\textnormal{sign}(x)+\frac12$. The bound of Theorem \ref{thm:robust_construction} becomes weaker and weaker. This contrasts with the bound of Theorem \ref{thm:construction}, which is completely independent of the pot function but requires initial ratings be exactly equal.
\end{remark}
}

\section{Upper bound for general $n$}
\label{sec:opt}
In this section, we wish to show the algorithm presented in Section \ref{sec:opt_strat_large_n} is nearly optimal.
Our strategy is to show that achieving a rating of $r$ requires many games to be played.
We start with a relaxation of the setup: instead of considering a discrete sequence of games resulting in a discrete sequence of player ratings,
we consider a continuous path in the space of possible player ratings.
Call a path $\rvec:[0,c]\to\R^n$ \textbf{valid} if there exists a finite sequence $0=t_0<\cdots<t_k=c$ such that for each $j\in[k]$ there exists $w,\ell\in[n]$ with
\[
\rvec'(t)=\e_w-\e_\ell\quad\forall\,t\in(t_{j-1},t_{j})
\]
where $\e_w,\e_\ell$ are elementary basis vectors.
For convenience, we write $\rvec'(t_{j-1})=\e_w-\e_\ell$ to make $\rvec'$ right-continuous.
In other words, $\rvec$ has constant speed and consists of line segments (called \textbf{edges}) along which only two coordinates change.
The points $\rvec(t_j)$ are referred to as \textbf{vertices}.
We will use the notation
\[A\to B\to C\]
to denote the path with vertices $A,B,C$.
A sequence of games corresponds naturally to a valid path:
let $\rvec(t_j)$ be the list of the players' ratings after game $j$.
The indices `$w$' and `$\ell$' conveniently correspond to the `winner' and `loser' for each game. In the event of a draw, the initially lower rated player is considered the winner as their rating increases. Note 
\begin{equation}
\label{eq:t_spacing}
    t_{j+1}-t_j = \text{amount of Elo transferred for game }j+1.
\end{equation}
Define a weight function
\begin{equation}
\label{eq:weight_function}
W(\mathbf x,\mathbf y)=\begin{cases}
    \frac1{\sigma(2-\inr{\mathbf x}{\mathbf y})} & \inr{\mathbf x}{\mathbf y}>0\\
    0&\inr{\mathbf x}{\mathbf y}\le0
\end{cases}
\end{equation}
and define the length of the path to be
\[
\len(\rvec)=
\int_0^c
W\pare{\rvec'(t),\rvec(t)}\wrt t.
\]
Whereas the Euclidean length of a path is simply $c\,\sqrt 2$, the value of $\len(\rvec)$ depends highly on the motion of the path.
When a valid path comes from a sequence of games, Lemma \ref{lem:length_to_k} shows $\len(\rvec)$ cannot be larger than $k$. The intuition is
that $\inr{\rvec'(t)}{\rvec(t)}=\rvec_w(t)-\rvec_\ell(t)$ records the difference in two players' ratings, so the Euclidean length of the $(j+1)$th edge is is $\sqrt2\sigma(-\inr{\rvec'(t_j)}{\rvec(t_j)})$, which approximately cancels with $w(\rvec'(t),\rvec(t))$
resulting in constant contribution.
\begin{lemma}
\label{lem:length_to_k}
Let $\rvec$ be a valid path arising from a sequence of $k$ games. Then
\[\len(\rvec)\le k
.\]
\end{lemma}
\begin{proof}
It suffices to show for each edge $\rvec(t_j)\to\rvec(t_{j+1})$ that
\[\len(\rvec(t_j)\to\rvec(t_{j+1}))\le1
\]
since summing over all edges would produce the final result.
Fix any $j$. Let $w,\ell$ be the indices such that
\[\rvec'(t)={\e_w-\e_\ell}\quad\forall\, t\in[t_{j},t_{j+1}).\]
Let $z=\inr{\rvec'(t_j)}{\rvec(t_j)}$
so that by (\ref{eq:t_spacing}) we have
\[\sigma(-z)=t_{j+1}-t_j.\]
Also note $\inr{\rvec'(t)}{\rvec(t)}=z+2(t-t_j)\le z+2\sig{-z}\le z+2$.
Then
\begin{align*}
    \len(\rvec(t_j)\to\rvec(t_{j+1}))
&=\int_{t_j}^{t_{j+1}}W(\rvec'(t),\rvec(t))\wrt t
\\&\le\int_{t_j}^{t_{j+1}}\frac1{\sig{2-\inr{\rvec'(t)}{\rvec(t)}}}\wrt t
\\&\le\int_{t_{j}}^{t_{j+1}}\frac1{\sig{-z}}\wrt t
\\&=1.
\end{align*}
\end{proof}
Lemma \ref{lem:length_to_k} means that instead of bounding $k$, we can bound $\len(\rvec)$. Directly doing so for any valid path is difficult, so we introduce two additional notions that allow us to make some additional restrictions on $\rvec$.
Call a valid path \textbf{ordered} if
\[\rvec_1(t)\ge\cdots\ge\rvec_n(t)\quad\forall\,t\in[0,c].\]
Each valid path corresponds to an ordered valid path of the same length by doing the following: each time $\rvec$ intersects a hyperplane of the form $x_w=x_\ell$, reflect the remaining part of the path across it.
This makes it possible to refer unambiguously to the $p$th highest rated player. It also means $\inr{\rvec'}\rvec$ does not change sign on each edge.
The last notion we introduce is that of an \textbf{upset}. This is when a lower rated player beats or draws with a higher rated player. Stated in terms of valid paths, an upset is an edge for which
\[\inr{\rvec'}{\rvec}=\rvec_w-\rvec_\ell<0\]
at the starting vertex.
We essentially show an optimal strategy, when converted to a ordered valid path, does not make use of upsets. The precise statement is given in Lemma \ref{lem:upset_free}. The intuition here is that upsets bring the players' ratings closer together, whereas in order to make one rating large you need the ratings to be spread out. The proof strategy is to take any valid path with upsets and convert it to one achieving a higher maximum rating without upsets. This lemma is the longest and most technical as it involves some casework in describing this conversion procedure.
\begin{lemma}
\label{lem:upset_free}
For each valid ordered path $\rvec:[0,c]\to\R^n$ there exists a valid ordered path $\tilde\rvec:[0,\tilde c]\to\R^n$ with
\[\len(\tilde\rvec)\le\len(\rvec)\quad\text{and}\quad
\max\tilde\rvec(\tilde c)\ge\max\rvec(c)
\] such that $\rvec'$ has no upsets.
\end{lemma}
\begin{proof}
Given $\rvec$, we start by constructing $\tilde\rvec$ with a slightly different property than the requirement of the theorem. We require $\rvec(c)=\tilde\rvec(\tilde c)$ and that all the upsets in $\tilde\rvec$ occur at the end. That is, once $\inr{\tilde\rvec'}{\tilde\rvec}$ becomes negative, it stays negative. Once we have done that, observe that upsets only decrease the value of $\max\tilde\rvec$ so we obtain a larger maximum value by truncating the path just before the upsets.

Fix vector $\mathbf u=(1,2,\cdots,n)$. Then each valid path $\rvec$ has an associated sequence $S(\rvec)$ of integers $S(\rvec)_j=\inr{\mathbf u}{\rvec'(t_j)}$ with one integer for each edge. Since that path is ordered, each integer is positive if and only if the corresponding edge is an upset. Equip the set of possible values of $S(\rvec)$ with the lexicographic ordering. Our method for constructing $\tilde\rvec$ is iterative, where each iteration strictly decreases $S(\rvec)$ and does not increase $\len(\rvec)$. The strict decreasing of $S(\rvec)$ guarantees that this procedure terminates after at most $n^{\#\text{ edges in }\rvec}$ steps.


We now describe one iteration. Given a path $\rvec$, locate three consecutive vertices $\rvec(t_j),\rvec(t_{j+1}),\rvec(t_{j+2})$ with
\[S(\rvec)_j>0\quad\text{and}\quad S(\rvec)_{j+1}<0.\]
That is, the two connecting edges are a upset followed by a non-upset.
We modify the path to obtain $\tilde\rvec$ by either deleting the middle vertex $\rvec(t_{j+1})$ from the path and adding an edge directly from $\rvec(t_j)$ to $\rvec(t_{j+2})$, or replacing the middle vertex with new vertex $\svec$. If $\rvec(t_{j+1})$ is deleted, $S(\rvec)$ is shortened and therefore decreased. If $\rvec(t_{j+1})$ is replaced by $\svec$, note that only the $j$th and $(j+1)$th entries in $S(\rvec)$ are affected. So requiring that
\[S(\tilde\rvec)_j<S(\rvec)_j\]
suffices for $S(\rvec)$ to decrease overall. Note the right hand side is positive, so it suffices for $S(\tilde\rvec)_j$ to be negative. 
We additionally require
\[\len(\tilde\rvec)\le\len(\rvec).\]
The modification differs based on how many coordinates change on those two edges. Edge edge modifies two coordinates, so this number can be 2, 3, or 4.

\textbf{Case of 4 coordinates}: this corresponds to the players of the two games being disjoint pairs. Intuitively, the order of the games is irrelevant to the outcome. Set $\svec=\rvec(t_j)-\rvec(t_{j+1})+\rvec(t_{j+2})$. Then $S(\tilde\rvec)_j=S(\rvec)_{j+1}<0$ and $\len(\rvec)=\len(\tilde\rvec)$.
\textbf{Case of 2 coordinates}: this corresponds to the same two players playing in consecutive games. Let those players be $w,\ell$ with $w<\ell$. Let $\Pi$ be the projection onto the $w,\ell$ coordinates. Then for some $a,b>0$,
\[
\Pi\rvec(t_j)=\bmat{x\\y}
\quad
\Pi\rvec(t_{j+1})=\bmat{x-a\\y+a}
\quad
\Pi\rvec(t_{j+2})=\bmat{x-a+b\\y+a-b}.
\]
Our modification simply deletes $\rvec(t_{j+1})$, so $S$ certainly decreases.
The original length is
\begin{align*}
\len\pare{\bmat{x\\y}\to\bmat{x-a\\y+a}\to\bmat{x-a+b\\y+a-b}}
  &=\int_0^b\frac1{\sig{2+  (y+a-t)-(x-a+t)  }}\wrt t
\\&=\int_0^b\frac1{\sig{2+  y-x+2a-2t  }}\wrt t.
\\&=\int_{-a}^{b-a}\frac1{\sig{2+  y-x-2t  }}\wrt t.
\end{align*}
If $a\ge b$, the new length is 0. If $a<b$, then
\begin{align*}
    \len\pare{\bmat{x\\y}\to\bmat{x-a+b\\y+a-b}}
  &=\int_0^{b-a}\frac1{\sig{2+  (y-t)-(x+t)  }}\wrt t
\\&=\int_0^{b-a}\frac1{\sig{2+  y-x-2t  }}\wrt t
\end{align*}
which is strictly less than the original length since the integrand is the same with a smaller range.
\textbf{Case of 3 coordinates}: this corresponds to one player playing two different opponents. Let $\Pi$ be the projection onto those three coordinates in ranked order, i.e.
\[\Pi\rvec(t_j)=\bmat{x\\y\\z}\]
for $x\ge y\ge z$.
There are six possible sign patterns of $\Pi\pare{\rvec(t_{j+2})-\rvec(t_j)}$, each with a different selection of $\svec$.

\textbf{Subcase $(+,-,+)$:}
For some $a,c>0$ we have \[\Pi\rvec(t_{j+2})=\bmat{x+a\\y-a-c\\z+c}.\] There are three possible values of $\Pi\rvec(t_{j+1})$ corresponding to the three possible upsets that can occur. No matter which one we observe, the selection of their replacement $\svec$ is the same.
\[
\Pi\rvec(t_{j+1})=
\bmat{x\\y-c\\z+c},
\bmat{x\\y-a-c\\z+a+c},
\bmat{x-c\\y\\z+c},
\quad
\svec=\bmat{x+a\\y-a\\z}
\]
The possible original lengths are as follows.
\begin{align*}L_1=\len\pare{\bmat{x\\y\\z}\to\bmat{x\\y-c\\z+c}\to\bmat{x+a\\y-a-c\\z+c}}
  &=\int_0^a \frac1{\sig{2+  (y-c-t)-(x+t)  }}\wrt t
\\&=\int_0^a \frac1{\sig{2+  y-x-c-2t  }}\wrt t,
\\
L_2=\len\pare{\bmat{x\\y\\z}\to\bmat{x\\y-a-c\\z+a+c}\to\bmat{x+a\\y-a-c\\z+c}}
  &=\int_0^a \frac1{\sig{2+  (z+a+c-t)-(x+t)  }}\wrt t
\\&=\int_0^a \frac1{\sig{2+  z-x+a+c-2t  }}\wrt t,
\\
L_3=\len\pare{\bmat{x\\y\\z}\to\bmat{x-c\\y\\z+c}\to\bmat{x+a\\y-a-c\\z+c}}
  &=\int_0^{a+c} \frac1{\sig{2+  (y-t)-(x-c+t)  }}\wrt t
\\&=\int_0^{a+c} \frac1{\sig{2+  y-x+c-2t  }}\wrt t
\\&=\int_{-c}^{a} \frac1{\sig{2+  y-x-c-2t  }}\wrt t
\end{align*}
The new length is
\[L_4=\len\pare{\bmat{x\\y\\z}\to\bmat{x+a\\y-a\\z}\to\bmat{x+a\\y-a-c\\z+c}}
=\int_0^a\frac1{\sig{2+  (y-t)-(x+t)  }}\wrt t
=\int_0^a\frac1{\sig{2+  y-x-2t  }}\wrt t.
\]
By monotonicity of $\sigma$, we automatically have $L_3\ge L_1\ge L_4$. We assume $\rvec$ is an ordered path, so in the $L_2$ case we have \(z+a+c\le y-a-c\) giving $L_2\ge L_4$ showing our modification did not increase the length. Also see that $\rvec(t_j)\to\svec$ is not an upset so $S(\tilde\rvec)_j<0$.

\textbf{Subcase $(-,+,-)$:}
For some $a,c>0$ we have \[\Pi\rvec(t_{j+2})=\bmat{x-a\\y+a+c\\z-c}.\] Again there are three possible values of $\Pi\rvec(t_{j+1})$ corresponding to the three possible upsets that can occur, and our selection of their replacement $\svec$ is the same regardless.
\[
\Pi\rvec(t_{j+1})=
\bmat{x-a\\y+a\\z},
\bmat{x-a-c\\y+a+c\\z},
\bmat{x-a\\y\\z+a},\quad
\svec=\bmat{x\\y+c\\z-c}.
\]
The original lengths are as follows.
\begin{align*}
L_1=\len\pare{\bmat{x\\y\\z}\to\bmat{x-a\\y+a\\z}\to\bmat{x-a\\y+a+c\\z-c}}
  &=\int_0^c \frac1{\sig{2+  (z-t)-(y+a+t)  }}\wrt t
\\&=\int_0^c \frac1{\sig{2+  z-y-a-2t  }}\wrt t,
\\
L_2=\len\pare{\bmat{x\\y\\z}\to\bmat{x-a-c\\y+a+c\\z}\to\bmat{x-a\\y+a+c\\z-c}}
  &=\int_0^c \frac1{\sig{2+  (z-t)-(x-a-c+t)  }}\wrt t
\\&=\int_0^c \frac1{\sig{2+  z-x+a+c-2t  }}\wrt t,
\\
L_3=\len\pare{\bmat{x\\y\\z}\to\bmat{x-a\\y\\z+a}\to\bmat{x-a\\y+a+c\\z-c}}
  &=\int_0^{a+c} \frac1{\sig{2+  (z-t)-(y+t)  }}\wrt t
\\&=\int_0^{a+c} \frac1{\sig{2+  z-y-2t  }}\wrt t,
\\&=\int_{-a}^{c} \frac1{\sig{2+  z-y-2a-2t  }}\wrt t.
\end{align*}
The new length is
\[L_4=\len\pare{\bmat{x\\y\\z}\to\bmat{x\\y+c\\z-c}\to\bmat{x-a\\y+a+c\\z-c}}
=\int_0^c\frac1{\sig{2+  (z-t)-(y+t)  }}\wrt t
=\int_0^c\frac1{\sig{2+  z-y-2t  }}\wrt t.
\]
Again by monotonicity of $\sigma$, we automatically have $L_3\ge L_1\ge L_4$. We assume $\rvec$ is an ordered path, so in the $L_2$ case we have \(y+a+c\le x-a-c\) giving $L_2\ge L_4$ as required. Also see that $\rvec(t_j)\to\svec$ is not an upset so $S(\tilde\rvec)_j<0$.

\textbf{Subcase $(+,+,-)$:}
For some $a,b>0$ we have
\[\Pi\rvec(t_{j+2})=\bmat{x+a\\y+b\\z-a-b}.\]
There is only one possible value of $\Pi\rvec(t_2)$ corresponding to an upset.
\[\Pi\rvec(t_{j+1})=\bmat{x-b\\y+b\\z},\quad\svec\bmat{x\\y+b\\z-b}.\]
The original length is
\begin{align*}
    L_1=\len\pare{\bmat{x\\y\\z}\to\bmat{x-b\\y+b\\z}\to\bmat{x+a\\y+b\\z-a-b}}
  &=\int_0^{a+b}\frac1{\sig{2+  (z-t)-(x-b+t)  }}\wrt t
\\&=\int_0^{a+b}\frac1{\sig{2+  z-x+b-2t  }}\wrt t
\\&=\int_0^b\frac1{\sig{2+  z-x+b-2t  }}\wrt t+\int_0^a\frac1{\sig{2+  z-x-b-2t  }}\wrt t
\end{align*}
and the new length is
\begin{align*}
L_2=\len\pare{\bmat{x\\y\\z}\to\bmat{x\\y+b\\z-b}\to\bmat{x+a\\y+b\\z-a-b}}
  &=\int_0^b\frac1{\sig{2+  (z-t)-(y+t)  }}\wrt t+\int_0^a\frac1{\sig{2+  (z-b-t)-(x+t)  }}\wrt t
\\&=\int_0^b\frac1{\sig{2+  z-y-2t  }}\wrt t+\int_0^a\frac1{\sig{2+  z-x-b-2t  }}\wrt t.
\end{align*}
Again since we assume the path is ordered, we have $x\ge y+b$ so $L_2\le L_1$ by monotonicity of $\sigma$. Also see that $\rvec(t_j)\to\svec$ is not an upset so $S(\tilde\rvec)_j<0$.

\textbf{Subcase $(+,-,-)$:}
For some $b,c>0$, we have
\[\Pi\rvec(t_{j+2})=\bmat{x+b+c\\y-b\\z-c}.\]
There is only one possible value of $\Pi\rvec(t_2)$ corresponding to an upset.
\[\Pi\rvec(t_{j+1})=\bmat{x\\y-b\\z+b},\quad\svec=\bmat{x+b\\y-b\\z}.\]
The original length is
\begin{align*}
L_1
   =\len\pare{\bmat{x\\y\\z}\to\bmat{x\\y-b\\z+b}\to\bmat{x+b+c\\y-b\\z-c}}
  &=\int_0^{b+c}\frac1{\sig{2+  (z+b-t)-(x+t)  }}\wrt t
\\&=\int_0^{b+c}\frac1{\sig{2+  z-x+b-2t  }}\wrt t
\\&
=\int_0^b\frac1{\sig{2+  z-x+b-2t  }}\wrt t
+\int_0^c\frac1{\sig{2+  z-x-b-2t  }}\wrt t
\end{align*}
and the new length is
\begin{align*}
L_2=\len\pare{\bmat{x\\y\\z}\to\bmat{x+b\\y-b\\z}\to\bmat{x+b+c\\y-b\\z-c}}
  &
=\int_0^b\frac1{\sig{2+  (y-t)-(x+t)  }}\wrt t
+\int_0^c\frac1{\sig{2+  (z-t)-(x+b+t)  }}\wrt t
\\&
=\int_0^b\frac1{\sig{2+  y-x-2t  }}\wrt t
+\int_0^c\frac1{\sig{2+  z-x-b-2t  }}\wrt t.
\end{align*}
Since we assume the path is ordered, we have $y-b\ge z$ so $L_2\le L_1$ once again. Also see that $\rvec(t_j)\to\svec$ is not an upset so $S(\tilde\rvec)_j<0$.

\textbf{Subcase $(-,+,+),(-,-,+)$:} For some $a,c$ we have
\[\Pi\rvec(t_{j+2})=\bmat{x-a\\y+a-c\\z+c}.
\]
Depending on the sign of $a-c$, there are two possible values of $\Pi\rvec(t_{j+1})$. In either case, the selection of $\svec$ is the same.
\[\Pi\rvec(t_{j+1})=\bmat{x-a\\y\\z+a},\bmat{x-c\\y\\z+c}\quad\svec=\bmat{x-a\\y+a\\z}.\]
In this case, both new edges are upsets so the new length is 0, but we cannot conclude $S(\tilde\rvec)_j<0$ as before. Instead, let $w_1<w_2<w_3$ be the indices of the three relevant coordinates and note \[S(\rvec)_j=w_1-w_3>S(\tilde\rvec)_j=w_1-w_2.\]
Finally note the sign patterns of $(+,+,+)$ and $(-,-,-)$ are not possible since the sum of all ratings is invariant.
Since in all possible cases we have a decrease of $S$ without an increase of length, this process terminates in a path that isn't longer and doesn't have any upsets followed by a non-upset. We end by truncating off any upsets at the end of the path.
\end{proof}
Lemma \ref{lem:upset_free} implies we can restrict our attention to valid ordered upset-free paths. The length of these paths can be bounded in terms of the following potential function.
\begin{definition}
    Define $\Phi:\R^n\to\R$ by
\[
\Phi(\svec)
=\magn{\svec}^2+\sum_{p=1}^{n-1}f(-2+\svec_p-\svec_{p+1})
\]
where $f$ is as in (\ref{eq:f_declr}).
\end{definition}
We seek to show $\Phi$ grows slowly as games are played, so that $\Phi$ of the end point of the path is upper bounded by (a multiple of) the length of the path.
Note that the first term $\magn \svec^2$ in this potential function is exactly the $\sum r^2$ expression appearing in the proofs of Theorems \ref{thm:construction} and \ref{thm:robust_construction}.
The intuitive idea behind the second term is that when player $p$ beats player $p+1$, the ratings $\svec_p$ and $\svec_{p+1}$ move $2\sigma(\svec_{p+1}-\svec_p)$ apart, which by the definition of $f$ means the corresponding $f$ term increases by just a constant amount. Lemma \ref{lem:phi_upper} makes this intuition precise.
\begin{lemma}
\label{lem:phi_upper}
For a valid ordered upset-free path $\rvec:[0,c]\to\R^n$ one has
\[\Phi(\rvec(c))-\Phi(\rvec(0))
\le\pare{2+2\sup_{z}z\,\sigma(2-z)}\cdot\len(\rvec)\]
\end{lemma}
\begin{proof}
By the fundamental theorem of calculus, it suffices to show that $$\deriv t{}\Phi(\rvec(t))\le\pare{2+2\sup_z z\sig{2-z}}W(\rvec'(t),\rvec(t)).$$
Since $\rvec$ is upset-free we have, $\inr{\rvec'(t)}{\rvec}\ge0$ so by the definition of $W$ (\ref{eq:weight_function}),
\[
W(\rvec'(t),\rvec(t))=\frac1{\sig{2-\inr{\rvec'(t)}{\rvec(t)}}}.
\]
The chain rule gives \begin{equation}\label{eq:chain_rule}\deriv t{}\Phi(\rvec(t))=\inr{\rvec'(t)}{\nabla \Phi(\rvec(t))}.\end{equation}
Since $\rvec$ is ordered in addition to being upset-free, $\rvec'(t)$ is always of the form $\rvec'(t)=\e_j-\e_{j+a}$ for positive $a$. That is, the winning player has a lower index than the losing player. 
The $j$th entry of the gradient can be computed directly as
\begin{align*}
\frac{\partial}{\partial\svec_j}\Phi(\svec)
  &=2\svec_j+f'(-2+ \svec_j-\svec_{j+1})\cdot\indicate{j<n}-f'(-2+ \svec_{j-1}-\svec_{j})\cdot\indicate{j>1}
\\&=2\rvec_j
+\frac{\indicate{j<n}}{\sig{2+ \svec_{j+1}-\svec_j}}
-\frac{\indicate{j>1}}{\sig{2+ \svec_{j}-\svec_{j-1}}},
\end{align*}
and in particular
\begin{align*}
\inr{\e_j-\e_{j+a}}{\nabla\Phi(\svec)}=
\pare{
\frac{\partial}{\partial\svec_j}\Phi(\svec)
-
\frac{\partial}{\partial\svec_{j+a}}\Phi(\svec)
}
  &\le2(\svec_j-\svec_{j+a})
+\frac{\indicate{j<n}}{\sig{2+ \svec_{j+1}-\svec_j}}
+\frac{\indicate{j+a>1}}{\sig{2+ \svec_{j+a}-\svec_{j+a-1}}}
\\&\le2(\rvec_j-\rvec_{j+a})
+\frac{\indicate{j<n}}{\sig{2+ \svec_{j+a}-\svec_j}}
+\frac{\indicate{j+a>1}}{\sig{2+ \svec_{j+a}-\svec_{j}}}
\\&\le
 \frac{2\sup_{z}z\,\sigma(2-z)}{\sig{2+ \svec_{j+a}-\svec_j}}
+\frac2{\sig{2+ \svec_{j+a}-\svec_j}}
\\&=
 \pare{2+2\sup_{z}z\,\sigma(2-z)}\cdot\frac1{\sig{2-\inr{\e_j-\e_{j+a}}{\svec}}}.
\\&=
 \pare{2+2\sup_{z}z\,\sigma(2-z)}\cdot w(\e_j-\e_{j+a},\svec).
\end{align*}
For $\svec=\rvec(t)$ and $\e_j-\e_{j+a}=\rvec'(t)$, this is exactly what we needed to show.
\end{proof}
The last piece of the puzzle is showing $\Phi(\svec)$ has to be large whenever one entry of $\svec$ is large.
Rapid growth of $f$ means if any two consecutive players have a large rating difference, $\Phi$ will be large. On the other hand, if consecutive players are close in rating, many entries of $\svec$ are close to its largest entry forcing $\magn\svec^2$ to be large.
\begin{lemma}
    \label{lem:phi_lower}
Let $\svec$ be any vector such that $\max\svec=\svec_1$ and at least one entry is negative. Then
\[\Phi(\svec)\ge\svec_1^2\quad\quad\text{and}\quad\quad\Phi(\svec)\ge\frac{\svec_1^3/8}{f^{-1}(\svec_1^2/4)+2}.\]
Furthermore, 
\[\Phi(\svec)\ge nf\pare{-2+\frac{\svec_1}{2n}}.\]
\end{lemma}
\begin{proof}
The first result is immediate from $\Phi(\svec)\ge\magn\svec^2\ge\svec_1^2$.
By the assumption on $\svec$, there exists $m$ such that
\[\svec_m\ge\svec_1/2\ge\svec_{m+1}.\]
Then
\begin{equation}
\label{eq:phi_part_1_lb}
\magn\svec^2\ge \frac m4\svec_1^2.
\end{equation}
By convexity of $f$,
\begin{align}
\label{eq:phi_part_2_lb}
\sum_{j=1}^{m}f\pare{-2+ \svec_j-\svec_{j+1}}
   \ge
\sum_{j=1}^{m}f\pare{-2+ \frac{\svec_1-\svec_{m+1}}m}
   \ge
mf\pare{-2+ \frac{\svec_1}{2m}}.
\end{align}
We claim (\ref{eq:phi_part_2_lb}) is monotonically decreasing in $m$.
To see this, note $f'$ is itself positive and monotonically increasing by monotonicity of $\sigma$. Then for $x>-2$,
\begin{align*}
f(x)=\int_0^xf'(t)\wrt t
\le xf'(x)
\le (x+2)f'(x)
\implies
0\le\frac{(x+2)f'(x)-f(x)}{(x+2)^2}=\deriv x{}\frac{f(x)}{x+2}.
\end{align*}
Setting $x=-2+\svec_1/2m$ and noting $x$ is itself monotonically decreasing in $m$ establishes the claim. $m$ is the index of a player, so we must have $m\le n$. This immediately implies by montonicity of (\ref{eq:phi_part_2_lb}) that
\[\Phi(\svec)\ge nf\pare{-2+\frac{\svec_1}{2n}}.\]
Combining (\ref{eq:phi_part_1_lb}) and (\ref{eq:phi_part_2_lb}) yields
\begin{align*}
\label{eq:phi_lb}
\Phi(\svec)
  &\ge\inf_m\pare{\frac m4\svec_1^2+mf\pare{-2+ \frac{\svec_1}{2m}}}
\\&\ge\inf_m\max\pare{\frac{m}4\svec_1^2,\,mf\pare{-2+ \frac{\svec_1}{2m}}}.
\end{align*}
The minimum of the maximum of two functions occurs when they intersect, which in this case is guaranteed to happen exactly once since $f(-2+\svec_1/2m)$ is monotone in $m$ and its range contains $[0,\infty)\ni\svec_1^2/4$. Therefore the minimizing $m$ is
\[m=\frac{\svec_1/2}{f^{-1}\pare{\svec_1^2/4}+2}.\]
The last result of the lemma follows immediately by plugging that $m$ into (\ref{eq:phi_part_1_lb}).
\end{proof}

Assembling the above lemmas together results in Theorem \ref{thm:opt}.
\begin{theorem}
\label{thm:opt}
Let $r$ be the highest rating achieved by a group of any number of players who play a total of $k$ games.
Suppose the pot function $\sigma$ satisfies
\[
C_1=\sup_{z}z\,\sigma(2-z)<\infty.
\]
Then
\[r\le 2n\cdot\pare{f^{-1}\pare{C_2\cdot{k}/{n}}+2}.\]
and
\[r\le
2\cdot C_2^{1/3}\cdot k^{1/3}\cdot\pare{f^{-1}\pare{C_2/4\cdot k}+2}^{1/3}
\]
for
\(C_2=2+2C_1+f(-2)\)
where $f$ is defined in (\ref{eq:f_declr}).
\end{theorem}
\begin{proof}
We may take $n\le k-1$ players without loss of generality
by simply ignoring players who who aren't connected to the player achieving the highest rating via some sequence of games.
Let $\rvec:[0,c]\to\R^n$ be the valid path corresponding to the sequence of games achieving rating $r=\max\rvec(c)$.
Perform all necessary reflections to make it ordered.
Lemma $\ref{lem:length_to_k}$ gives
\begin{equation}
\label{finaleq:len_k}
    \len(\rvec)\le k.
\end{equation}
Lemma \ref{lem:upset_free} constructs $\tilde\rvec$ without any upsets satisfying
\begin{equation}
\label{finaleq:r_r}
    r=\max\rvec(c)\le\max\tilde\rvec(\tilde c)
\end{equation}
and
\begin{equation}
\label{finaleq:len_len}
    \len(\tilde\rvec)\le\len(\rvec).
\end{equation}
Since $\tilde\rvec$ is ordered and upset-free, Lemma \ref{lem:phi_upper} implies
\[\Phi(\tilde\rvec(\tilde c))-\Phi(\tilde\rvec(0))\le(2+2C_1)\,\len(\tilde\rvec).\]
Note $\Phi(\tilde\rvec(0))=\Phi(0)=(n-1)f(-2)$ so
\begin{equation}
\label{finaleq:phi_len}
\Phi(\tilde\rvec(\tilde c))\le(2+2C_1)\,\len(\tilde\rvec)+(n-1)f(-2).
\end{equation}
Set $C_2=2+2C_1+f(-2)$ and assemble the chain of inequalities:
\begin{align*}
\Phi(\tilde\rvec(\tilde c))
  &\underset{(\ref{finaleq:phi_len})}\le(2+2C_1)\len(\tilde\rvec)+(n-1)f(-2)
\\&\underset{(\ref{finaleq:len_len})}\le(2+2C_1)\len(\rvec)+(n-1)f(-2)
\\&\underset{(\ref{finaleq:len_k})}\le(2+2C_1)k+(n-1)f(-2)
\\&\underset{}\le (2+2C_1+f(-2))\cdot k
\\&\underset{}= C_2\cdot k
\stepcounter{equation}\tag{\theequation}\label{finaleq:chain}
\end{align*}
The three lower bounds on $\Phi$ provided by Lemma \ref{lem:phi_lower} are each used in different ways. First, the $n$-dependent bound implies
\[
nf\pare{-2+\frac{r}{2n}}
\underset{(\ref{finaleq:r_r})}\le
nf\pare{-2+\frac{\max\tilde\rvec(\tilde c)}{2n}}
\underset{\ref{lem:phi_lower}}\le
\Phi(\tilde\rvec(\tilde c))
\underset{(\ref{finaleq:chain})}\le C_2\cdot k.
\]
Since $r$ appears only once in the equation, we can simply rearrange to solve for $r$. In particular,
\[r\le 2n\pare{f^{-1}\pare{C_2\cdot\frac{k}{n}}+2}.\]
giving the first result of the theorem.
The second $n$-independent bound
gives
\begin{equation}
\label{finaleq:thing_phi}
\frac{r^3/8}{f^{-1}(r^2/4)+2}
\underset{(\ref{finaleq:r_r})}\le
\frac{(\max\tilde\rvec(\tilde c))^3/8}{f^{-1}((\max\tilde\rvec(\tilde c))^2/4)+2}
\underset{\ref{lem:phi_lower}}\le
\Phi(\tilde\rvec(\tilde c))
\underset{(\ref{finaleq:chain})}\le C_2\cdot k.
\end{equation}
This bounds $k$ in terms of a function of $r$.
The last application of Lemma \ref{lem:phi_lower} converts this into a bound on $r$ in terms of $k$.
\[r^2
\underset{(\ref{finaleq:r_r})}\le
\pare{\max\tilde\rvec(\tilde c)}^2
\underset{\ref{lem:phi_lower}}\le
\Phi(\tilde\rvec(\tilde c))
\underset{(\ref{finaleq:chain})}\le C_2\cdot k.
\]
This can be plugged back into the left-hand side of (\ref{finaleq:thing_phi}),\[
C_2\cdot k\ge
\frac{r^3/8}{f^{-1}\pare{r^2/4}+2}
\ge
\frac{r^3/8}{f^{-1}\pare{C_2/4\cdot k}+2}.
\]
Rearranging gives
\[
r^3\le 8\cdot C_2\cdot k\cdot\pare{f^{-1}\pare{C_2/4\cdot k}+2}.\]
Taking cube roots establishes the final result.
\end{proof}

\begin{remark}
    The probabilistic interpretation (\ref{eq:prob_interp}) of the Elo system described in the introduction lends itself naturally to the requirement that $z\,\sigma(2-z)$ be bounded. In particular, take $\sig{z}=Pr(\eta<x)+\frac12\Pr(\eta=x)$ to be the symmetrized cumulative distribution function of a symmetric random variable $\eta$ \textbf{with finite expectation}.
Note
\[\sup_z z\,\sigma(2-z)=\sup_z\,(2+z)\,\sigma(-z)\le2+\sup_zz\,\sigma(-z).\]
It's clear that the supremum on the right will be achieved for $z\ge0$.
Since $\eta$ has finite expectation, we may apply Markov's inequality to obtain
\[
z\sig{-z}\le z\Pr(\eta\le-z)=\frac12\cdot z\cdot \Pr(\abs\eta\ge z)\le\frac12\E\abs\eta<\infty.
\]
This shows $\E\abs\eta<\infty$ is sufficient. However it isn't strictly necessary; for instance, $z\,\sigma(2-z)$ just barely does not diverge for a Cauchy random variable. However, one \textit{does} need finite $(1-\eps)$th moment: suppose $\sig{-z}\le c/z$. Then
\[
\E\abs{\eta}^{1-\eps}
=\int\Pr\pare{\abs{\eta}^{1-\eps}\ge t}\wrt t
=2\int\sig{-t^{1/(1-\eps)}}\wrt t
\le2\int{t^{-1/(1-\eps)}}\wrt t<\infty\]
for each $\eps\in(0,1)$.
\end{remark}
\section{Discussion}
\label{sec:dis}
\subsection{Remaining questions}
\label{sec:remaining}
Can one close the gap between the upper and lower bounds in the $n=\omega(1)$ regime?
This could occur by finding a better strategy than the ones in Section \ref{sec:opt_strat_large_n}, for instance by slowly increasing $\delta$ as games are played, or by tightening the analysis in Theorem \ref{thm:opt}.
It's possible that for heavy-tailed $\sigma$ the lower bound is too loose, but for light-tail $\sigma$ the upper bound is.
Why is there a jump in the upper bound from $f(k)^{1/3}$ to $f(k/n)$ at $n=k^{1/3}$? Can one find a bound that smoothly crosses the phase transition?
In Theorem \ref{thm:robust_construction}, one cannot totally specify the initial set of ratings, but has to allow the possibility that the initial ratings are all flipped. This appears very strongly to be an artifact of the analysis; can one prove a version that allows you to assign any initial ratings subject only to the constraint that the average rating is non-negative?

The highest Elo problem also bares a passing resemblance to the Toda lattice: players are particles whose positions on the real line is given by their ratings; they exhibit a repulsive force when a higher rated player beats a lower rated player and an attractive force if vice versa. Is it possible to use tools from solid state physics to analyze this problem?

Another interesting variant is to constrain the number of \textit{rounds} of games. In a tournament, many games will be happening parallel. Each round consists of any number of games, subject only to the constraint that each player only participates in at most one game each round. \textbf{Given $n$ players and $k$ \textit{rounds}, what's the highest someone may be rated after all games have finished?}

\subsection{Connection to maximum overhang}

One may notice the jump from $\log k$ to $k^{1/3}$ as one increases $n$ for $\sigma(-z)=\frac1{1+e^{-z}}$. This may be reminiscent of the maximum overhang problem.
In that problem, one places $k$ unit-length bricks on top of each other at the edge of a table, attempting to achieve the largest overhang possible. The classic solution if one is restricted to placing only a single brick at each height (the `single-wide' setting) achieves an overhang of $\Theta(\log k)$ units. However, if one is allowed to place as many bricks at each height as one likes (the `multi-wide' setting), the optimal solution achieves an overhang of $\Theta(k^{1/3})$ \cite{overhang2}. This connection does not appear to be a coincidence! The proof of optimality in the maximum overhang problem uses a reduction to a `mass movement problem.'
But in fact, the highest Elo problem for the particular $\sigma(z)=\max(0,\min(1,z/2+1/2))$ can be reduced to nearly the same problem!

The mass movement problem described by \cite{overhang2} is as follows:
consider some finite number of piles of mass placed on the real line. A valid `move' takes some unit interval and rearranges the mass within that interval so that the center of mass is unchanged. Negative mass is allowed. Formally, a `signed distribution' $\mu$ is a finite linear combination of dirac delta functions $\delta$. A valid `move' replaces $\mu$ with $\mu+\nu$ where $\nu$ is itself a signed distribution whose support is contained in some unit interval and $\int\nu=0$.
A sequence of moves $\nu_1\cdots,\nu_k$ corresponds naturally to a sequence of signed distributions $\mu_0,\cdots,\mu_k$.
They place an additional `weight-constraint' on allowed sequences. In particular, the number of $j$ for which $\nu_j$ is allowed to have support to the right of any threshold $T$ is at most $\max_{0\le i\le k}\int_T^\infty\mu_i$.
The goal is to use at most $k$ moves to move the distribution $k\cdot\delta$ to a distribution with a unit of mass as far right as possible.

To reduce highest Elo to mass movement, let $r_p(j)$ be the rating of player $p$ after $j$ games and set
\[\mu_j(x)=\sum_p\delta(x-r_p(j)).\]
Then consider a game where $p$ beats $q$.
Set $z=r_p-r_q$. By the selection of $\sigma(z)=\max(0,\min(z/2+1/2))$, this only produces a change in ratings if $z\le1$. Then $\mu$ changes by performing the following move:
\begin{align*}
    \nu(x)
=\delta(x-r_p-\sig{-z})
+\delta(x-r_q+\sig{-z})
-\delta(x-r_p)
-\delta(x-r_q).
\end{align*}
Suppose $z\ge -1$. Then the support of $\nu$ is contained in $[r_q-\sig{-z},r_p+\sig{-z}]$, which is an interval of length $z+2\sig{-z}=1$ as required.
Furthermore, $\nu$ only has support above $T$ when there's a player above $T$ Elo guaranteeing the weight-constraint.
The objective is starting with $n\cdot\delta$, produce via $k$ moves a distribution with a unit of mass (i.e. at least one player) as far right as possible (i.e. with the highest rating possible). Recall that for the purposes of an upper bound, we make take $n=k$, establishing the correct objective.

There is, however, one big discrepancy: the case of $z<-1$. In particular, in the highest Elo problem it is possible for a player rated $-10$ to beat a player rated $10$, garnering a full point of Elo. This does not correspond to a valid move, since the support of $\nu$ would be $\{-10,-9,9,10\}$. Intuitively, such moves should not ultimately help push mass far away, and indeed in Lemma \ref{lem:upset_free} we show for a different relaxation of highest Elo that they don't help.

A second discrepancy is that mass movement places no restriction on the amount of mass moved in each move, whereas for highest Elo we always have $\magn{\nu}_1=2$. This would seem to imply an upper bound of the `$k$ rounds' variant mentioned at the end of Section \ref{sec:remaining}. Can that argument be made to work? Conversely, is there a variant or generalization of maximum overhang that can be analyzed using the proof of Theorem \ref{thm:opt}?

\pagebreak
\bibliographystyle{plain} 
\bibliography{refs} 

\end{document}